\documentclass[10pt]{article}
\usepackage{amsfonts}
\usepackage{amsmath}
\usepackage{mathrsfs}
\usepackage{mathrsfs,amscd,amssymb,amsthm,amsmath,psfrag,bm,graphicx}

\makeatletter

\renewcommand{\@seccntformat}[1]{{\csname the#1\endcsname}{\normalsize .}\hspace{.5em}}
\makeatother

\def \[{\begin{equation}}
\def \]{\end{equation}}

\newtheorem{thm}{Theorem}[section]
\newtheorem{prop}{Proposition}

\newtheorem{claim}{Claim}
\newtheorem{fact}{Fact}
\newtheorem{lem}[thm]{Lemma}

\newenvironment{kst}
{\setlength{\leftmargini}{2\parindent}
 \begin{itemize}
 \setlength{\itemsep}{-1.1mm}}
{\end{itemize}}

\newenvironment{wst}
{\setlength{\leftmargini}{1.5\parindent}
 \begin{itemize}
 \setlength{\itemsep}{-1.1mm}}
{\end{itemize}}

\begin{document}
\setlength{\baselineskip}{15pt}
\begin{center}{\Large \bf On the spectral moment of graphs with given clique number\footnote{Financially supported by the National Natural
Science Foundation of China (Grant Nos. 11071096, 11271149) and the Special Fund for Basic Scientific Research of Central Colleges (CCNU11A02015).}}

\vspace{4mm}

{\large Shuna Hu,\ \ Shuchao Li\footnote{E-mail: 374289328@qq.com (S.N. Hu),\ lscmath@mail.ccnu.edu.cn (S.C.
Li)},\ \ Xixi Zhang}\vspace{2mm}

Faculty of Mathematics and Statistics,  Central China Normal
University, Wuhan 430079, P.R. China
\end{center}
\noindent {\bf Abstract}: \ Let $\mathscr{L}_{n,t}$ be the set of all $n$-vertex connected graphs with clique number $t$\,($2\leq t\leq n)$. For $n$-vertex connected graphs with given clique number, lexicographic ordering
by spectral moments ($S$-order) is discussed in this paper. The first $\sum_{i=1}^{\lfloor\frac{n-t-1}{3}\rfloor}(n-t-3i)+1$ graphs with $3\le t\le n-4$, and the last few graphs, in the $S$-order, among $\mathscr{L}_{n,t}$ are characterized. In addition, all graphs in $\mathscr{L}_{n,n}\bigcup\mathscr{L}_{n,n-1}$ have an $S$-order; for the cases $t=n-2$ and $t=n-3$ the first three and the first seven graphs in the set $\mathscr{L}_{n,t}$ are characterized, respectively.

\vspace{2mm} \noindent{\it Keywords}: Spectral moment;  Clique number;  Chromatic number

\vspace{2mm}

\section{\normalsize Introduction}\setcounter{equation}{0}
All graphs considered here are finite, simple and connected. Undefined terminology and notation may be referred to \cite{D-I}.
Let $G=(V_G,E_G)$ be a simple undirected graph with $n$ vertices. $G-v$, $G-uv$ denote the graph obtained from $G$ by deleting vertex $v \in V_G$, or edge
$uv \in E_G$, respectively (this notation is naturally extended if more than one vertex or edge is deleted). Similarly,
$G+uv$ is obtained from $G$ by adding an edge $uv \not\in E_G$. For $v\in V_G,$ let
$N_G(v)$ (or $N(v)$ for short) denote the set of all the adjacent vertices of $v$ in $G$ and $d_G(v)=|N_G(v)|$. A \textit{pendant vertex} of $G$ is a vertex of degree 1.

Let $A(G)$ be the adjacency matrix of a graph $G$ with $\lambda_1(G),\lambda_2(G),\dots,\lambda_n(G)$ being its eigenvalues in non-increasing order. The number $\sum_{i=1}^n\lambda_i^k(G)\, (k=0,1,\dots,n-1)$ is called the \textit{$k$th spectral moment} of $G$, denoted by $S_k(G)$. Let $S(G)=(S_0(G), S_1(G), \ldots, S_{n-1}(G))$ be the sequence of spectral moments of $G$. For two graphs $G_1, G_2$, we shall write $G_1=_sG_2$ if $S_i(G_1)=S_i(G_2)$ for $i=0,1,\dots,n-1$. Similarly, we have $G_1\prec_sG_2\, (G_1$ comes before $G_2$  in the $S$-order) if for some $k\, (1\leq k\leq {n-1})$, we have  $S_i(G_1)=S_i(G_2)\, (i=0,1,\dots,k-1)$  and $S_k(G_1)<S_k(G_2)$. We shall also write $G_1\preceq_sG_2$ if $G_1\prec_sG_2$ or $G_1=_sG_2$.  $S$-order has been used in producing graph catalogs (see \cite{C-G}), and for a more general setting of spectral moments one may be referred to \cite{C-R-S1}.

Recently, investigation on $S$-order of graphs has received increasing attention.
For example, Cvetkovi\'c and Rowlinson \cite{C-R-S3} studied
the $S$-order of trees and unicyclic graphs and characterized the first and the last graphs, in the $S$-order,
of all trees and all unicyclic graph with given girth, respectively. Wu and Fan \cite{D-F} determined the first
and the last graphs, in the $S$-order, of all unicyclic graphs and bicyclic graphs, respectively. Pan et al. \cite{X-F-P}
gave the first $\sum_{k=1}^{\lfloor\frac{n-1}{3}\rfloor}(\lfloor\frac{n-k-1}{2}\rfloor-k+1)$ graphs apart from an $n$-vertex path, in the $S$-order, of all trees with $n$ vertices. Wu and Liu \cite{C-M-M} determined the last $\lfloor\frac{d}{2}\rfloor+1$ graphs, in the $S$-order,
among all $n$-vertex trees of diameter $d\, (4 \le d \le n-3)$. Pan et al. \cite{B-Z1} identified the
last and the second last graphs, in the $S$-order, of quasi-trees. Cheng, Liu and Liu identified the last
$d + \lfloor \frac{d}{2}\rfloor-2$ graphs, in the $S$-order, among all $n$-vertex unicyclic graphs of diameter $d$. Cheng and Liu \cite{C-M-M1} determined the
last few graphs, in the $S$-order, among all trees with $n$ vertices and $k$ pendant vertices. Li and Song \cite{Li-S} identified the last $n$-vertex tree with a given degree sequence in the $S$-order. Consequently, the last trees in the $S$-order among the sets of all trees of order $n$ with the largest degree, the leaves number, the independence number and the matching number was also determined, respectively. Li and Zhang \cite{Li-Z} determined the first, the second, the last and the second last graphs in the $S$-order among the set of all graphs with given number of cut edges.

On the other hand, there are many Tur\'an-type extremal problems, i.e., given a forbidden graph $F$,
determine the maximal number of edges in a graph on $n$ vertices that does not contain a copy of
$F$. It states that among $n$-vertex graphs not containing a clique of size $t + 1$, the complete $t$-partite
graph $T_{n,t}$ with (almost) equal parts, which is called Tur\'an graph, has the maximum number of edges.
Spectral graph theory has similar Tur\'an extremal problems which determine the largest (or smallest)
eigenvalue of a graph not containing a subgraph $F$.  Nikiforov explicitly
proposed to study the general Tur\'an problems in [12-16]. For example, he
\cite{13} determined the maximum spectral radius of graphs without paths of given length and presented a comprehensive survey on these topics; see \cite{15}. In addition, Sudakov et al. \cite{17} presented a generalization of Tur\'an Theorem in terms of Laplacian eigenvalues, whereas He et al. \cite{D-M} gave a generalization of Tur\'an Theorem in terms of signless Laplacian eigenvalues.

Motivated by Tur\'an-type extremal problems, we investigate in this paper
the spectral moments of $n$-vertex graphs with given clique number,
which may be regarded as a part of spectral extremal theory.
For $2\le t\le n$, let $\mathscr{L}_{n,t}$ be the set of all $n$-vertex connected graphs with clique number $t$. We give the first $\sum_{i=1}^{\lfloor\frac{n-t-1}{3}\rfloor}(n-t-3i)+1$ graphs with $3\le t\le n-4$, and the last few graphs, in the $S$-order, among $\mathscr{L}_{n,t}$. In addition, all graphs in $\mathscr{L}_{n,n}\bigcup\mathscr{L}_{n,n-1}$ have an $S$-order; for the cases $t=n-2$ and $t=n-3$ the first three and the first seven graphs in the set $\mathscr{L}_{n,t}$ are characterized respectively. We prove these results in Section 3. According to the relationship between the clique number and the chromatic number of graphs, we study the $S$-order of graphs with given chromatic number in Section 4. In Section 2, we give some preliminaries which are useful for the proofs of our
main results.

\section{\normalsize Preliminaries}\setcounter{equation}{0}
Throughout we denote by $P_n,\,S_n,\, C_n$ and $K_n$ a path, a star, a cycle and a complete graph on $n$ vertices, respectively. An $F$-\textit{subgraph} of $G$ is a subgraph of $G$, which is isomorphic to the graph $F$. Let $\phi_G(F)$ (or $\phi(F)$ for short) be the number of all $F$-subgraphs of $G$. The notation
$G\nsupseteq F$ means that $G$ does not contain $F$ as its subgraph.

Further on we will need the following lemmas.
\begin{lem}[\cite{B-Z1}]
The $k$th spectral moment of $G$ is equal to the number of closed walks of length $k$.
\end{lem}

Let $H_1, H_2, \ldots, H_{23}$ be the graphs as depicted in Fig. 1, which will be used in Lemma 2.2.
\begin{figure}[h!]
\begin{center}
\psfrag{a}{$H_{10}$} \psfrag{b}{$H_{11}$}
\psfrag{c}{$H_{12}$} \psfrag{d}{$H_{13}$}
\psfrag{e}{$H_{14}$} \psfrag{f}{$H_{15}$}
\psfrag{h}{$H_{16}$} \psfrag{1}{$H_1$}
\psfrag{i}{$H_{17}$} \psfrag{j}{$H_{18}$}
\psfrag{k}{$H_{19}$} \psfrag{m}{$H_{20}$}
\psfrag{m}{$H_{21}$} \psfrag{n}{$H_{22}$}
\psfrag{o}{$H_{23}$}\psfrag{l}{$H_{20}$}
\psfrag{2}{$H_2$} \psfrag{3}{$H_3$}
\psfrag{4}{$H_4$} \psfrag{5}{$H_5$}
\psfrag{6}{$H_6$} \psfrag{7}{$H_7$}
\psfrag{8}{$H_8$} \psfrag{9}{$H_9$}
\includegraphics[width=130mm]{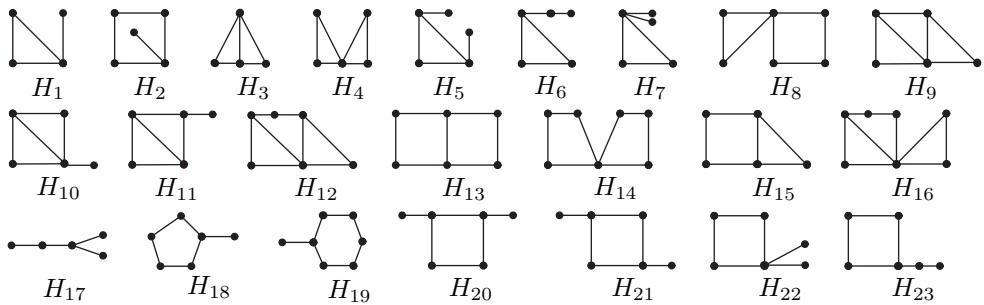}\\
\caption{Graphs $H_1, H_2, \ldots, H_{22}$ and $H_{23}$.}
\end{center}
\end{figure}
\begin{lem}[\cite{C-R-S2}]
Given a connected graph $G$, $S_0(G)=n, S_1(G)=l, S_2(G)=2m, S_3(G)=6t$, where $n, l, m, t$ denote the number of vertices, the number of loops, the number of edges and the number of triangles contained in $G$, respectively.
\end{lem}
\begin{lem}
For every graph $G$, we have
\begin{wst}
\item[{\rm (i)}] $S_4(G )=2 \phi(P_2)+4 \phi(P_3)+8 \phi(C_4)$ {\rm (\cite{C-R-S1})}.

\item[{\rm (ii)}]$S_5(G )=30 \phi(C_3)+10 \phi(H_1)+10 \phi(C_5)$ {\rm (\cite{C-R-S1})}.

\item[{\rm (iii)}]$S_6(G )=2 \phi(P_2)+12 \phi(P_3)+6 \phi(P_4)+12 \phi(K_{1,3})+12 \phi(H_2)+ 36 \phi(H_3)+ 24 \phi(H_4)+ 24\phi(C_3)+ 48 \phi(C_4)+ 12 \phi(C_6)$ {\rm (\cite{C-M-M})}.

\item[{\rm (iv)}]
    $S_7(G)=126\phi(C_3)+84\phi(H_1)+28\phi(H_7)+14\phi(H_5)+14\phi(H_6)+112\phi(H_3)+42\phi(H_{15})+28\phi(H_8)+70\phi(C_5)+14\phi(H_{18})+14\phi(C_7).$

\item[{\rm (v)}]
    $S_8(G)=2 \phi(P_2)+28 \phi(P_3)+32 \phi(P_4)+8 \phi(P_5)+72 \phi(K_{1,3})+16 \phi(H_{17})+48 \phi(K_{1,4})+168 \phi(C_3)+64 \phi(H_1)+464 \phi(H_3)+384 \phi(H_4) +96 \phi(H_{15})+96 \phi(H_{10})+48 \phi(H_{11})+80 \phi(H_{12})+32 \phi(H_{16})+264 \phi(C_4)+24 \phi(H_9)+112 \phi(H_2)+16 \phi(H_{23})+16 \phi(H_{20})+16 \phi(H_{21})+32 \phi(H_{22})+32 \phi(H_{13})+32 \phi(H_{14})+528 \phi(K_4)+96 \phi(C_6)+16 \phi(H_{19})+16 \phi(C_8).$
\end{wst}
\end{lem}
\begin{proof}
(iv)\ By Lemma 2.1 and we note that vertices that belong to a closed walk of length $7$  induce in $G$ a subgraph isomorphic to $C_3,\,H_1,\, H_7,\,H_5,\,H_6,\,H_3,\,H_{15},\,H_8,\,C_5,\,H_{18},\,C_7$. By using matlab, we can obtain  the number of closed walks of length 7 which span these subgraphs is 126, 84, 28, 14, 14, 112, 42, 28, 70, 14 and 14, respectively,  then (iv) follows.

(v)\ By Lemma 2.1 and we note that vertices that belong to a closed walk of length $8$ induce in $G$ a subgraph isomorphic to $P_2,\,P_3,\, P_4,\,P_5,\,K_{1,3},\,H_{17},\,K_{1,4},\,C_3,\,H_1,\,H_3,\,H_4,\,H_{15},\,H_{10},\,H_{11},\,H_{12},\,H_{16},\,C_4,\, H_9,\linebreak H_2,\, H_{23},\,H_{20}, H_{21},\, H_{22},\,H_{13},\,H_{14},\,K_4,\,C_6,\,H_{19},\,C_8$. By using matlab, we can obtain  the number of closed walks of length 8 which span these subgraphs is 2, 28, 32, 8, 72, 16, 48, 168, 64, 464, 384, 96, 96, 48, 80, 32, 264, 24, 112, 16, 16, 16, 32, 32, 32, 528, 96, 16 and 16, respectively, then (v) follows immediately.
\end{proof}

A connected subgraph $H$ of $G$ is called a \textit{tree-subgraph} (or \textit{cycle-subgraph}) if $H$ is a tree (or contains
at least one cycle). Let $H$ be a proper subgraph of $G$, we call $H$ an \textit{effective graph} for $S_k(G)$ if $H$ contains a closed walk of length $k$. Set $\mathscr{T}_k(G)=\{T:\, T$ is a tree-subgraph of $G$ with $|E_T|\leq\frac{k}{2}\}$; $\mathscr{T}'_k(G)=\{W:\, W$ is a cycle-subgraph of $G$ with $|E_W|\leq k$\}; $\mathscr{A}_k(G)=\{T:\, T$ is a tree-subgraph of $G$ and and it is an effective graph for $S_k(G)\}$; $\mathscr{A}'_k(G)=\{W:\, W$ is a cycle-subgraph of $G$ and it is an effective graph for $S_k(G)\}.$ It is easy to see that $\mathscr{A}_k(G) \cap \mathscr{A}'_k(G)=\emptyset$. By Lemma 2.1 we have
\begin{prop}
Given a graph $G$,  the set of all effective graphs for $S_k(G)$ is $\mathscr{A}_k(G) \cup \mathscr{A}'_k(G).$ In particular, if $k$ is odd, then $\mathscr{A}_k(G )=\emptyset.$
\end{prop}
\begin{lem} [\cite{C-M-M}]
Let $G$ be a non-trivial connected graph with $u\in V_G$. Suppose that two paths of lengths $a$,\,$b\,(a\geq b\geq 1)$ are attached to $G$ by their end vertices at $u$, respectively, to form $G^{*}_{a,b}$. Then $G^{*}_{a+1,b-1}\prec_sG^{*}_{a,b}$.
\end{lem}

Let $G$ and $H$ be two graphs with $u\in V_G$ and $v\in V_H$, we shall denote by $Gu\cdot vH$ the graph obtained from $G$ and $H$ by identifying $u$ and $v$.
\begin{lem}[\cite{C-M-M}]\label{lem2.5}
Let $G,\,H$ be two non-trivial connected graphs with $u,\,v\in V_G$ and $w\in V_H$. If $d_G(u) <d_G(v)$, then $Gu\cdot wH \prec_s Gv\cdot wH$.
\end{lem}

\section{\normalsize On the $S$-order among $\mathscr{L}_{n,t}$}\setcounter{equation}{0}

In this section, we study the $S$-order among $\mathscr{L}_{n,t} \,(2\leq t\leq n).$ In view of Lemma 2.2 the first few graphs in the $S$-order among $\mathscr{L}_{n,2}$ must be $n$-vertex trees. Fortunately,  on the other hand, Pan et al. \cite{X-F-P} identified the first $\sum_{k=1}^{\lfloor\frac{n-1}{3}\rfloor}\left(\lfloor\frac{n-k-1}{2}\rfloor-k+1\right)+1$ graphs, in the $S$-order, of all trees with $n$ vertices; these
$\sum_{k=1}^{\lfloor\frac{n-1}{3}\rfloor}\left(\lfloor\frac{n-k-1}{2}\rfloor-k+1\right)+1$ trees are also the first $\sum_{k=1}^{\lfloor\frac{n-1}{3}\rfloor}\left(\lfloor\frac{n-k-1}{2}\rfloor-k+1\right)+1$ graphs in the $S$-order among  $\mathscr{L}_{n,2}$. We will not repeat it here.

A graph $G\not\supseteq F$ on $n$ vertices with the largest possible number of edges is called \textit{extremal} for $n$ and $H$; its number of edges is denoted by ex$(n,F)$.  The following theorem tells us the Tur\'an graph $T_{n,t}$ is indeed extremal for $n$ and $K_{t+1},$ and as such unique.
\begin{thm} [Tur\'an 1941]
For all integers $t,n$ with $t>1,$ every graph $G\not\supseteq K_{t+1}$ with $n$ vertices and {\rm ex}$(n,K_{t+1})$ edges is a $T_{n,t}.$
\end{thm}
\begin{thm}
For all integers $t,n$ with $t>1,$ every graph $G\in \mathscr{L}_{n,t} \backslash \{T_{n,t}\}$, one has $G\prec_s T_{n,t}$.
\end{thm}
\begin{proof}
Note that, for graph $G \in \mathscr{L}_{n,t}\setminus \{T_{n,t}\}$, one has $S_i(G)=S_i(T_{n,t}),\, i=0,\,1$. By Lemma 2.2 and Theorem 3.1, we have $S_2(G)=2|E_G|<2|E_{T_{n,t}}|=S_2(T_{n,t})$. Hence, $G\prec_sT_{n,t}.$
\end{proof}

Now assume that $(V_1,\,V_2,\,\ldots,\,V_t)$ is a partition of $T_{n,t}$ with $n=kt+r\,(0\leq r<t),$ where $|V_i|=k$ if $i=1,\,2,\,\ldots,\,t-r$ and $|V_i|=k+1$ otherwise. For $u,\,v\in V_{T_{n,t}}$, let

\medskip\noindent
$\bullet$  $T_{n,t}^1$ be the graph obtained by deleting the edge $uv$ from $T_{n,t}$, where $u\in V_i, v\in V_j$,\, $1\leq i\not=j\leq t-r;$

\medskip\noindent
$\bullet$  $T_{n,t}^2$ be the graph obtained by deleting the edge $uv$ from $T_{n,t}$, where $u\in V_i,v\in V_j$,\, $1\leq i\leq t-r$ and $t-r+1\leq j\leq t$;

\medskip\noindent
$\bullet$  $T_{n,t}^3$ be the graph obtained by deleting the edge $uv$ from $T_{n,t}$, where $u\in V_i, v\in V_j$,\, $t-r+1\leq i\not=j\leq t$.\vspace{2mm}

In particular, if $t=n-1$ then $|V_1|=|V_2|=\cdots =|V_{n-2}|=1$ and $|V_{n-1}|=2$. In this case, for convenience we assume that $V_i=\{v_i\}$ for $i=1,\,2,\,\ldots,\,n-2$ and $V_{n-1}=\{v_{n-1},u\}$. 
Let
$$
T_i:= T_{n,n-1}-\{uv_1,\,uv_2,\,\cdots,\,uv_{i-1},\, uv_i\},
$$
where $i=1,\,2,\,\ldots,\,n-3.$ It is straightforward to check that $\mathscr{L}_{n,n-1}=\{T_{n,n-1},\, T_1,\, T_2,\, \ldots,\, T_{n-3}\}.$

\begin{thm}Among the set of graphs $\mathscr{L}_{n,t}$ with $3\le t \le n-1.$
\begin{wst}
\item[{\rm (i)}]If $n=t+1$, then all graphs in the set $\mathscr{L}_{n,n-1}$ have the following $S$-order:
$T_{n-3}\prec_s T_{n-4}\prec_s \cdots \prec_s T_i \prec_s \cdots\prec_s T_2\prec_sT_1\prec_s T_{n,n-1}.$
\item[{\rm (ii)}]If $n=kt$ with $3\leq t\leq\frac{n}{2}$, then for all $G\in \mathscr{L}_{n,t}\backslash\{T_{n,t},\,T_{n,t}^1\}$ one has $G\prec_s T_{n,t}^1\prec_s T_{n,t}.$
\item[{\rm (iii)}]If $n=kt+1$ with $3\leq t\leq\frac{n}{2}$, then for all $G\in \mathscr{L}_{n,t}\backslash\{T_{n,t},\,T_{n,t}^1,\,T_{n,t}^2\}$ one has $G\prec_s T_{n,t}^1\prec_s T_{n,t}^2\prec_s T_{n,t}.$
\item[{\rm (iv)}]If $n=kt+r$ with $3\leq t\leq\frac{n}{2},\, r=t-1$ or $\frac{n+1}{2}\leq t\leq n-2$, then for all $G\in \mathscr{L}_{n,t}\backslash\{T_{n,t},\,T_{n,t}^2,\,T_{n,t}^3\}$ one has $G \prec_s T_{n,t}^2 \prec_s T_{n,t}^3\prec_s T_{n,t}.$
\item[{\rm (v)}]If $n=kt+r$ with $4\leq t\leq\frac{n}{2},\, 2\leq r\leq t-2$, then for all $G\in \mathscr{L}_{n,t}\backslash\{T_{n,t},\,T_{n,t}^1,\,T_{n,t}^2,\,T_{n,t}^3\}$ one has $G\prec_s T_{n,t}^1\prec_s T_{n,t}^2 \prec_s T_{n,t}^3\prec_s T_{n,t}.$
\end{wst}
\end{thm}
\begin{proof}
(i)\ It follows directly by Lemma 2.2.

(ii)\ For all $G\in \mathscr{L}_{n,t}\backslash\{T_{n,t},\,T_{n,t}^1\}$, by Lemma 2.2 $S_i(G)=S_i(T_{n,t}^1)=S_i(T_{n,t}^2)$ for $i=0,\,1$. Note that $n=kt$ with $3\leq t\leq\frac{n}{2}$, hence by the definition of $T_{n,t}^1$, it is the unique graph in $\mathscr{L}_{n,t}$ satisfying the number of its edges equals to $|E_{T_{n,t}}|-1,$ which implies that for all $G\in \mathscr{L}_{n,t}\backslash\{T_{n,t},\,T_{n,t}^1\}$, we have $|E_G|<|E_{T_{n,t}^1}|,$ i.e., $S_2(G)<S_2(T_{n,t}^1)<S_2(T_{n,t})$. Hence, (ii) holds.

(iii)\ For any $G\in \mathscr{L}_{n,t}\backslash\{T_{n,t},\,T_{n,t}^1,\,T_{n,t}^2\}$, one has $S_i(G)=S_i(T_{n,t}^1)=S_i(T_{n,t}^2)$
for $i=0,\,1$. By the definition of $T_{n,t}^1,\, T_{n,t}^2$, we know they are just the two graphs in $\mathscr{L}_{n,t}$ satisfying $|E_{T_{n,t}^1}|=|E_{T_{n,t}^2}|=|E_{T_{n,t}}|-1,$ which implies that for all $G\in \mathscr{L}_{n,t}\backslash\{T_{n,t},\,T_{n,t}^1,\,T_{n,t}^2\}$, we have $|E_G|<|E_{T_{n,t}^1}|=|E_{T_{n,t}^2}|,$ i.e., $S_2(G)<S_2(T_{n,t}^1)=S_2(T_{n,t}^2)<S_2(T_{n,t})$. In order to complete the proof, it suffices to show $S_3(T_{n,t}^1)<S_3(T_{n,t}^2)$. In fact,
$$
 S_3(T_{n,t}^1)-S_3(T_{n,t}^2)=6(\phi_{T_{n,t}^1}(C_3)-\phi_{T_{n,t}^2}(C_3))=-6<0.
$$
Hence, (iii) holds.

(iv) and (v) can be proved by a similar discussion as in the proof of (iii). We omit the procedure here.
\end{proof}

In the following, we are to determine the first few graphs, in the $S$-order, among $\mathscr{L}_{n,t}\,(3\leq t\leq n-1).$ Let $J_1, J_2, \ldots, J_{10}$ be the $n$-vertex graphs as depicted in Fig. 2.
\begin{figure}[h!]
\begin{center}
\psfrag{a}{$u$} \psfrag{k}{$K_{n-2}$} \psfrag{m}{$K_{n-3}$}
\psfrag{1}{$J_1$} \psfrag{2}{$J_2$} \psfrag{3}{$J_3$}
\psfrag{4}{$J_4$}\psfrag{5}{$J_5$} \psfrag{6}{$J_6$} \psfrag{7}{$J_7$}
\psfrag{8}{$J_8$}\psfrag{9}{$J_9$} \psfrag{i}{$J_{10}$}
\includegraphics[width=120mm] {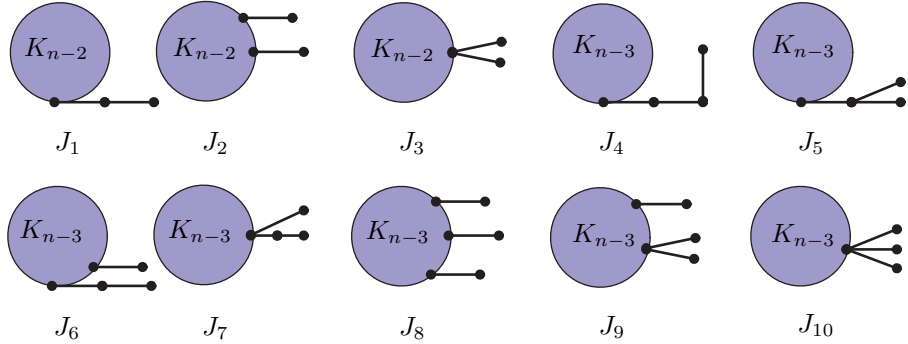}\\
\caption{$n$-vertex graphs $J_1, J_2,\ldots, J_{10}$.}
\end{center}
\end{figure}
\begin{thm}Among the set of all graphs $\mathscr{L}_{n,t}$.
\begin{kst}
\item[{\rm (i)}]If $t=n-2\ge 3,$ the first three graphs, in the $S$-order, among $\mathscr{L}_{n,n-2}$ are
$J_1\prec_s J_2 \prec_s J_3.$
\item[{\rm (ii)}]If $t=n-3\ge 3,$ the first seven graphs, in the $S$-order, among $\mathscr{L}_{n,n-3}$ are
$J_4\prec_s J_5 \prec_s J_6 \prec_s J_7 \prec_s J_8 \prec_s J_9 \prec_s J_{10}.$
\end{kst}
\end{thm}
\begin{proof}
(i)\ For all $G\in\mathscr{L}_{n,n-2}\setminus\{J_1, J_2, J_3\}$, one has not only $S_k(G)=S_k(J_1)=S_k(J_2)=S_k(J_3)$ for $k=0,1$
but also $|E_G|>|E_{J_1}|=|E_{J_2}|=|E_{J_3}|$, hence $S_2(J_1)=S_2(J_2)=S_2(J_3)<S_2(G)$, i.e., $J_i \prec_s G$ for $i=1,2,3$.
By Lemma 2.5, we have $J_1\prec_s J_2\prec_s J_3.$ Therefore, (i) holds.

(ii)\ For all $G\in\mathscr{L}_{n,n-3}\setminus\{J_4, J_5, J_6, J_7, J_8, J_9, J_{10}\}$, one has $S_k(G)=S_k(J_4)=\cdots =S_k(J_{10})$ for $k=0,1$
and $|E_G|>|E_{J_4}|=\cdots = |E_{J_{10}}|$. Hence $S_2(J_4)=S_2(J_5)=\cdots = S_2(J_{10})<S_2(G)$, i.e., $J_i \prec_s G$ for $i=4,5,\ldots, 10$. It is easy to see that $S_3(J_4)=S_3(J_5)=\cdots = S_3(J_{10}).$ Further on,
by Lemma 2.3(i)-(ii), we have
\begin{eqnarray}
  S_4(J_4) &=& 2|E_{J_4}|+8\phi_{J_4}(C_4)+4\left[(n-4){n-4\choose2}+{n-3\choose2}+2{2\choose2}\right],\label{eq:3.1}\\
  S_4(J_5) &=& 2|E_{J_5}|+8\phi_{J_5}(C_4)+4\left[(n-4){n-4\choose2}+{3\choose2}\right],\label{eq:3.2}\\
  S_4(J_6) &=& 2|E_{J_6}|+8\phi_{J_6}(C_4)+4\left[(n-5){n-4\choose2}+2{n-3\choose2}+{2\choose2}\right].\label{eq:3.3}
\end{eqnarray}

By (\ref{eq:3.1}) and (\ref{eq:3.2}), we have $S_4(J_4)-S_4(J_5)=-4<0$, hence $J_4\prec_s J_5.$ Similarly, we can show that
$J_6\prec_s J_7, J_8 \prec_s J_9 \prec_s J_{10}.$

By (\ref{eq:3.2}) and (\ref{eq:3.3}), we have $S_4(J_5)-S_4(J_6)=-4(n-6).$ If $n>6$, obviously we have
$S_4(J_5)-S_4(J_6)<0$, i.e., $J_5 \prec_s J_6$ for $n>6$. If $n=6,
S_4(J_5)-S_4(J_6)=0,$  whereas $S_5(J_5)-S_6(J_6)=10(\phi_{J_5}(U_4)-\phi_{J_6}(U_4))=-10<0$,
i.e., $J_5 \prec_s J_6$ for $n=6$. Hence, we obtain that $J_5 \prec_s J_6$ for $n\ge 6$. Similarly, we can
also show that $J_7 \prec_s J_8$ for $n\ge 6.$

Hence, we obtain $J_4\prec_s J_5 \prec_s J_6 \prec_s J_7 \prec_s J_8 \prec_s J_9 \prec_s J_{10} \prec_s G$ for $G\in \mathscr{L}_{n,n-3}\setminus\{J_4, J_5, J_6, J_7, J_8,\linebreak J_9, J_{10}\}$, as desired.
\end{proof}

\begin{figure}[h!]
\begin{center}
\psfrag{a}{$u$} \psfrag{b}{$K_t$}\psfrag{A}{$K_t^{n-t-i}$}
\psfrag{B}{$K_t^{n-t-i}v_j\cdot w_0P_{i+1}$}\psfrag{C}{$K_t^{n-t-i}u\cdot w_0P_{i+1}$}
\psfrag{0}{$v_0$} \psfrag{1}{$v_1$}
\psfrag{2}{$w_0$}\psfrag{3}{$w_1$}\psfrag{l}{$w_0$}
\psfrag{5}{$w_2$}\psfrag{t}{$w_i$}\psfrag{i}{$v_j$}
\psfrag{j}{$v_j$} \psfrag{4}{$v_{n-t-i}$}
\includegraphics[width=130mm]{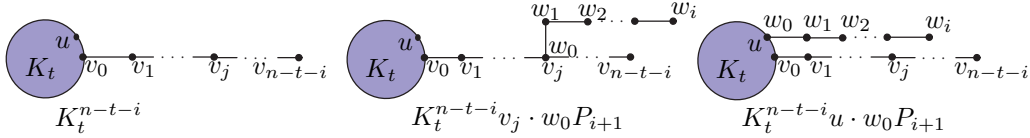}\\
\caption{Graphs $K_t^{n-t-i}, K_t^{n-t-i}v_j\cdot w_0P_{i+1}$ and $K_t^{n-t-i}u\cdot w_0P_{i+1}$ with some vertices labelled.}
\end{center}
\end{figure}

Given a path $P_{i+1}=w_0w_1w_2\ldots w_i$, set $G_j^i:=K_t^{n-t-i}v_j\cdot w_0P_{i+1}$ and $H^i:=K_t^{n-t-i}u\cdot w_0P_{i+1}$ (see Fig. 3), where $3\leq t\leq n-2,\,0\leq j\leq n-t-2i$ and $1\leq i\leq \lfloor\frac{n-t}{2}\rfloor$.
\begin{thm}
Among the set of all graphs $\mathscr{L}_{n,t}$ with $3\leq t\leq n-4,$ the first $\sum_{i=1}^{\lfloor\frac{n-t-1}{3}\rfloor}(n-t-3i)+1$
graphs, in the $S$-order, are $K_t^{n-t}\prec_sG^1_{n-t-2}\prec_sG^1_{n-t-3}\prec_s\cdots\prec_sG^1_3\prec_sG^1_2\prec_sG^2_{n-t-4}\prec_sG^2_{n-t-5}\prec_s\cdots\prec_sG^2_3\prec_s\cdots
\prec_s G^i_{n-t-2i} \prec_s G^i_{n-t-2i-1}\prec_s\cdots \prec_s G_{i+1}^i \prec_s G^i_{n-t-2(i+1)} \prec_s \cdots \prec_s
G^{\lfloor\frac{n-t-1}{3}\rfloor}_{\lfloor\frac{n-t-1}{3}\rfloor+1}.$
\end{thm}
\begin{proof}
Note that for any connected graph $G$, one has $G\prec_s G+e$, where $e\not\in E_G.$ Hence, for $3\leq t\leq n-4,$ the first graph in the $S$-order among $\mathscr{L}_{n,t}$ is obtained from $K_t$ by attaching some trees to some vertices of $K_t,$ in view of Lemma \ref{lem2.5}, the first few graphs in the $S$-order among $\mathscr{L}_{n,t}$ is just the kite graph $K_t^{n-t}.$ Furthermore, it suffices for us to consider the set of graphs $\mathscr{A}=\{G: G$ is an $n$-vertex graph obtained by attaching some trees to $K_t$ such that $G$ contains just two pendant vertices\}. It is easy to see that $\mathscr{A}$ is a subset of $\mathscr{L}_{n,t}$ and
$$
 \mathscr{A}=\left\{G_j^i:\,0\leq j\leq n-t-2i,\,1\leq i\leq \left\lfloor\frac{n-t}{2}\right\rfloor\right\}\bigcup\left\{H^i:\,1\leq i\leq \left\lfloor\frac{n-t}{2}\right\rfloor\right\}.
$$

We first show the following claims.
\begin{claim}
$G_j^i\prec_sH^{i'} \prec_s G_0^{i''}$, where $1\leq j\leq n-t-2i,\,1\leq i\leq \lfloor\frac{n-t-1}{2}\rfloor,\,1\leq i', i''\leq \lfloor\frac{n-t}{2}\rfloor;$
\end{claim}
\noindent{\bf Proof of Claim 1.}\ Note that, $1\leq j\leq n-t-2i,\,1\leq i\leq \lfloor\frac{n-t-1}{2}\rfloor,\,1\leq i',i''\leq \lfloor\frac{n-t}{2}\rfloor,$  $S_k(H^{i'})=S_k(G_j^i)=S_k(G_0^{i''})$ for $k=0,\,1,\,2,\,3.$
By Lemma 2.3(i), we have
\begin{eqnarray}
  S_4(H^{i'}) &=& 2m(H^{i'})+8\phi_{H^{i'}}(C_4)+4\left[(t-2){t-1\choose2}+2{t\choose2}+n-t-2\right],\label{eq:3.4}\\
  S_4(G_j^i) &=&2m(G_j^i)+8\phi_{G_j^i}(C_4)+4\left[(t-1){t-1\choose2}+{t\choose2}+{3\choose2}+n-t-3\right],\label{eq:3.5}\\
  S_4(G_0^{i''}) &=& 2m(G_0^{i''})+8\phi_{G_0^{i''}}(C_4)+4\left[(t-1){t-1\choose2}+{t+1\choose2}+n-t-2\right].\label{eq:3.6}
\end{eqnarray}
Note that $\phi_{H^{i'}}(C_4)=\phi_{G_j^i}(C_4)=\phi_{G_0^{i''}}(C_4)$,
hence (\ref{eq:3.4})-(\ref{eq:3.6}) gives
\[\label{eq:3.7}
 S_4(G_0^{i''})-S_4(H^{i'})=4,\, S_4(H^{i'})-S_4(G_j^i)=4(t-3).
\]
In view of (\ref{eq:3.7}), if $3< t\leq n-4$, we obtain
$G_j^i\prec_s H^{i'}\prec_s G_0^{i''}$; if $t=3$, then $S_4(H^{i'})=S_4(G_j^i)<S_4(G_0^{i''})$
and $S_5(H^{i'})-S_5(G_j^i)=10\phi_{H^{i'}}(U_4)-10\phi_{G_j^i}(U_4)=10>0,$ hence $G_j^i\prec_sH^{i'}\prec_s G_0^{i''}$ for $t=3$.
Hence $G_j^i\prec_s H^{i'}\prec_s G_0^{i''}$ for $3\leq t\leq n-4$, as desired. \qed

\begin{claim}
$G^i_{j+1}\prec_s G^i_j,$ where $1\leq i\leq\lfloor\frac{n-t-1}{2}\rfloor$, $1\leq j< n-t-2i.$
\end{claim}
\noindent{\bf Proof of Claim 2.}\ \
Note that, for $1\leq i\leq\lfloor\frac{n-t-1}{2}\rfloor$, $1\leq j<j+1\leq n-t-2i$, it is routine to check that $S_k(G_j^i)=S_k(G_{j+1}^i),\,$  $k\in \{0,\,1,\,2,\,3,\,4\}.$ In what follows we consider $k\ge 5$. On the one hand, for $5\leq k\leq 2j+3$, it is easy to see that for any $W\in\mathscr{A}'_k(G_j^i)$, there exists $W'\in\mathscr{A}'_k(G_{j+1}^i)$ such that $W\cong W'$, and vice versa. Hence,
\[\label{eq:3.8}
   \mathscr{A}'_k(G_j^i)=\mathscr{A}'_k(G_{j+1}^i),\ \ \ \  k=5,6,\ldots, 2j+3.
\]
In what follows, we distinguish our discussion in the following two possible cases.

\vspace{2mm}

{\bf Case 1.}\ \ $1\le j\leq\frac{n-t-i-1}{2}$. In this case, if $5\leq k\leq 2j+2$, then for any $T\in\mathscr{A}_k(G_j^i)$, there exists $T'\in \mathscr{A}_k(G_{j+1}^i)$ such that $T\cong T'$ for $5\leq k\leq 2j+2$, and vice versa. Note that if $k$ is odd, then $\mathscr{A}_k(G)=\emptyset$, hence $\mathscr{A}_{2j+3}(G_j^i)=\mathscr{A}_{2j+3}(G_{j+1}^i)=\emptyset$. Therefore,
\[\label{eq:3.9}
   \mathscr{A}_k(G_j^i)=\mathscr{A}_k(G_{j+1}^i),\ \ \ \  k=5,6,\ldots, 2j+3.
\]
By (\ref{eq:3.8}), (\ref{eq:3.9}) and Proposition 1, we obtain
\[\notag
S_k(G_j^i)=S_k(G_{j+1}^i), \ \ \ \ k=5,6,\ldots, 2j+3.
\]

If $1\le j<\frac{n-t-i-1}{2}$, note that $k=2j+4$, then for any $W\in\mathscr{A}'_{2j+4}(G_j^i)$, there exists $W'\in\mathscr{A}'_{2j+4}(G_{j+1}^i)$ such that $W\cong W'$, and vice versa. Hence, $\mathscr{A}'_{2j+4}(G_j^i)=\mathscr{A}'_{2j+4}(G_{j+1}^i).$ Notice that for any $T\in \mathscr{A}_{2j+4}(G_j^i),\, T'\in \mathscr{A}_{2j+4}(G_{j+1}^i)$, it is routine to check
\begin{itemize}
  \item if $|E_T\cap E_{K_t}| =0$, then $
\phi_{G_j^i}(T)-\phi_{G_{j+1}^i}(T')=\phi_{G_j^i}(P_{j+3})-\phi_{G_{j+1}^i}(P_{j+3})=-1;$
  \item if $|E_T\cap E_{K_t}| =1$, then $
\phi_{G_j^i}(T)-\phi_{G_{j+1}^i}(T')=\phi_{G_j^i}(P_{j+3})-\phi_{G_{j+1}^i}(P_{j+3})=t-1;$
  \item if $|E_T\cap E_{K_t}| \ge 2$, then $
\phi_{G_j^i}(T)-\phi_{G_{j+1}^i}(T')=0.$
\end{itemize}
Hence,
$$
|\mathscr{A}_{2j+4}(G_j^i)|-|\mathscr{A}_{2j+4}(G_{j+1}^i)|=\phi_{G_j^i}(T)-\phi_{G_{j+1}^i}(T')=\phi_{G_j^i}(P_{j+3})-\phi_{G_{j+1}^i}(P_{j+3})=t-2\geq1>0.
$$
By Proposition 1, we have
\begin{eqnarray*}
  S_{2j+4}(G_j^i)-S_{2j+4}(G_{j+1}^i)=(2j+4)(\phi_{G_j^i}(P_{j+3})-\phi_{G_{j+1}^i}(P_{j+3})) =  (2j+4)(t-2)>0,
\end{eqnarray*}
which implies that $G_{j+1}^i\prec_sG_j^i$ for $1\le j< \frac{n-t-i-1}{2}$.\vspace{2mm}

If $j=\frac{n-t-i-1}{2}$, by a similar was as above we can obtain that
\begin{eqnarray*}
   |\mathscr{A}'_{2j+4}(G_j^i)|-|\mathscr{A}'_{2j+4}(G_{j+1}^i)|&=&0, \\[3pt]
  |\mathscr{A}_{2j+4}(G_j^i)|-|\mathscr{A}_{2j+4}(G_{j+1}^i)|&=& (2j+4)(t-1)>0.
\end{eqnarray*}
By Proposition 1, we get $S_{2j+4}(G_j^i)>S_{2j+4}(G_{j+1}^i)$, i.e.,  $G^i_{j+1}\prec_sG^i_j$ holds for $j=\frac{n-t-i-1}{2}$.\vspace{2mm}

{\bf Case 2.}\  $\frac{n-t-i-1}{2}<j\leq n-t-2i$. \vspace{2mm}

Note that $j>\frac{n-t-i-1}{2}$, hence $2(n-t-i-j)+1<2j+3$. In view of (\ref{eq:3.8}), we have
$$
  \mathscr{A}_k'(G_j^i)=\mathscr{A}_k'(G_{j+1}^i),\ \ \ \  k=5,6,\ldots, 2(n-t-i-j)+1.
$$
Furthermore, for $5\leq k\leq 2(n-t-i-j)+1$, we have for all $T\in\mathscr{A}_k(G_j^i)$, there exists $T'\in \mathscr{A}_k(G_{j+1}^i)$ such that $T\cong T'$ and vice versa. Hence,
$$
  \mathscr{A}_k(G_j^i)=\mathscr{A}_k(G_{j+1}^i),\ \ \ \  k=5,6,\ldots, 2(n-t-i-j)+1.
$$
By Proposition 1, $S_k(G_j^i)=S_k(G_{j+1}^i)$ holds for $5\le k\le 2(n-t-i-j)+1$. For
$k=2(n-t-i-j+1)$, by a similar discussion as in the proof of the Case 1,  we can obtain that
\begin{eqnarray*}
   |\mathscr{A}'_{2(n-t-i-j+1)}(G_j^i)|-|\mathscr{A}'_{2(n-t-i-j+1)}(G_{j+1}^i)|&=&0, \\[3pt]
  |\mathscr{A}_{2(n-t-i-j+1)}(G_j^i)|-|\mathscr{A}_{2(n-t-i-j+1)}(G_{j+1}^i)|&=& \phi_{G_j^i}(T)-\phi_{G_{j+1}^i}(T')=1.
\end{eqnarray*}
By Proposition 1, $S_{2(n-t-i-j+1)}(G_j^i)-S_{2(n-t-i-j+1)}(G_{j+1}^i)=2(n-t-i-j+1)\cdot 1>0.$ Hence, $G^i_{j+1}\prec_sG^i_j$ holds for $\frac{n-t-i-1}{2}<j\leq n-t-2i$.

By Cases 1 and 2, Claims 2 holds. This completes the proof. \qed

\begin{claim}
$G^i_{i+1}\prec_sG^{i+1}_{n-t-2(i+1)},$ where $1\leq i<\lfloor\frac{n-t-1}{3}\rfloor$.
\end{claim}
\noindent{\bf Proof of Claim 3.}\  By a similar discussion as in the proof of Claims 1 and 2, we can show that $S_k(G^i_{i+1})=S_k(G^{i+1}_{n-t-2(i+1)})$ for $0\leq k\leq 2i+3$ and $S_{2i+4}(G^{i+1}_{n-t-2(i+1)})-S_{2i+4}(G^i_{i+1})=2i+4>0$. Hence, $G^i_{i+1}\prec_sG^{i+1}_{n-t-2(i+1)}$. \qed\vspace{2mm}

By Claims 1-3, we have the following fact.
\begin{fact}
The set $\mathscr{B}:=\{K_t^{n-t},\, G^1_{n-t-2},\, G^1_{n-t-3},\, \cdots, G^1_3,\, G^1_2,\, G^2_{n-t-4},\, G^2_{n-t-5},\, \cdots, G^2_3, \cdots, G^i_{n-t-2i}, \linebreak G^i_{n-t-2i-1}, \cdots, G_{i+1}^i,\,  G^i_{n-t-2(i+1)},\, \cdots G^{\lfloor\frac{n-t-1}{3}\rfloor}_{\lfloor\frac{n-t-1}{3}\rfloor+1}\}$ consists of $\sum_{i=1}^{\lfloor\frac{n-t-1}{3}\rfloor}(n-t-3i)+1$ graphs and they are in the following $S$-order:
$K_t^{n-t}\prec_sG^1_{n-t-2}\prec_sG^1_{n-t-3}\prec_s\cdots\prec_sG^1_3\prec_sG^1_2\prec_sG^2_{n-t-4}\prec_sG^2_{n-t-5}\prec_s\cdots\prec_s G^2_3\prec_s \cdots
\prec_s G^i_{n-t-2i} \prec_s G^i_{n-t-2i-1} \prec_s \cdots \prec_s G_{i+1}^i \prec_s G^i_{n-t-2(i+1)} \prec_s \cdots \prec_s
G^{\lfloor\frac{n-t-1}{3}\rfloor}_{\lfloor\frac{n-t-1}{3}\rfloor+1},$ where $3\leq t\leq n-4.$
\end{fact}

\begin{claim}Among $\mathscr{A}\setminus \mathscr{B},$ one has
$G^{\lfloor\frac{n-t}{3}\rfloor}_{\lfloor\frac{n-t}{3}\rfloor}\preceq_sG^i_j,$ where $1\leq j\leq n-t-2i,\, j\leq i\leq \lfloor\frac{n-t-1}{2}\rfloor.$
\end{claim}
\noindent{\bf Proof of Claim 4.}\  For a fixed $j$ in $\{1,\,2, \ldots, n-t-2i\}$,  there does not exist $G^i_j$ satisfying $i\ge j>\lfloor\frac{n-t}{3}\rfloor$. By Lemma 2.2, we know $G^{j}_j\prec_sG^i_j$ for all $ 1\leq j\leq \lfloor\frac{n-t}{3}\rfloor,\, j< i \leq \lfloor\frac{n-t-1}{2}\rfloor$. Hence, according to the $S$-order, the first graph in $\{G_j^i:\ 1\leq j\leq n-t-2i,\, j\leq i\leq \lfloor\frac{n-t-1}{2}\rfloor\}$ is just the first graph in $\{G^{j}_j:\ 1\leq j\leq \lfloor\frac{n-t}{3}\rfloor\}$. In what follows, we are to determine the first graph in the $S$-order among $\{G^{j}_j:\ 1\leq j\leq \lfloor\frac{n-t}{3}\rfloor\}$.

Note that $1\leq i+1\leq \lfloor\frac{n-t}{3}\rfloor-1$, hence $1\leq i+1\leq \frac{n-t-3}{3}.$ By a similar discussion as in the proof of Claim 1, we have $S_k(G_i^i)=S_k(G_{i+1}^{i+1})$ for $0\leq k\leq 2i+3$ and $S_{2i+4}(G_i^i)-S_{2i+4}(G_{i+1}^{i+1})= (2i+4)(t-3)$. Hence, if $t>3$ we obtain that  $S_{2i+4}(G_i^i)>S_{2i+4}(G_{i+1}^{i+1})$; if $t=3$, then $S_{2i+4}(G_i^i)=S_{2i+4}(G_{i+1}^{i+1})$. Furthermore, if $t=3$, $\mathscr{A}_{2i+5}(G_i^i)=\mathscr{A}_{2i+5}(G_{i+1}^{i+1})=\emptyset$ and for all $W\in \mathscr{A'}_{2i+5}(G_i^i),\, W'\in \mathscr{A'}_{2i+5}(G_{i+1}^{i+1})$, we have $\phi_{G_i^i}(W)-\phi_{G_{i+1}^{i+1}}(W')=1$. By Proposition 1, $S_{2i+5}(G_i^i)-S_{2i+5}(G_{i+1}^{i+1})\geq1>0$.
Hence, we obtain
$$
   G^{i+1}_{i+1}\prec_sG^i_i,\ \ \ \ 1\leq i< \lfloor\frac{n-t}{3}\rfloor,
$$
which implies that $G^{\lfloor\frac{n-t}{3}\rfloor}_{\lfloor\frac{n-t}{3}\rfloor}\preceq_sG^i_j$ for all $1\leq j\leq n-t-2i,\, j\leq i \leq\lfloor\frac{n-t-1}{2}\rfloor$, as desired. \qed

\begin{claim}
$G^{\lfloor\frac{n-t-1}{3}\rfloor} _{\lfloor\frac{n-t-1}{3}\rfloor+1}\prec_sG^{\lfloor\frac{n-t}{3}\rfloor} _{\lfloor\frac{n-t}{3}\rfloor}.$
\end{claim}
\noindent{\bf Proof of Claim 5.}\  Note that
$\lfloor\frac{n-t-1}{3}\rfloor=\lfloor\frac{n-t}{3}\rfloor$ or
$\lfloor\frac{n-t-1}{3}\rfloor=\lfloor\frac{n-t}{3}\rfloor-1$, for
the latter case, $\lfloor\frac{n-t}{3}\rfloor=\frac{n-t}{3}$. Hence,
if $\lfloor\frac{n-t-1}{3}\rfloor=\lfloor\frac{n-t}{3}\rfloor$, then
let $i=j=\lfloor\frac{n-t-1}{3}\rfloor$. By Claim 2,
$G^{\lfloor\frac{n-t-1}{3}\rfloor}
_{\lfloor\frac{n-t-1}{3}\rfloor+1}\prec_sG^{\lfloor\frac{n-t-1}{3}\rfloor}
_{\lfloor\frac{n-t-1}{3}\rfloor}$, i.e.,
$G^{\lfloor\frac{n-t-1}{3}\rfloor}
_{\lfloor\frac{n-t-1}{3}\rfloor+1}\prec_sG^{\lfloor\frac{n-t}{3}\rfloor}
_{\lfloor\frac{n-t}{3}\rfloor}$. If
$\lfloor\frac{n-t-1}{3}\rfloor=\frac{n-t}{3}-1$, by Lemma 2.4,
$G^{\frac{n-t}{3}-1}_{\frac{n-t}{3}}\prec_s
G^{\frac{n-t}{3}}_{\frac{n-t}{3}},$ i.e.,
$G^{\lfloor\frac{n-t-1}{3}\rfloor}
_{\lfloor\frac{n-t-1}{3}\rfloor+1}\prec_sG^{\lfloor\frac{n-t}{3}\rfloor}
_{\lfloor\frac{n-t}{3}\rfloor},$ as desired. \qed\vspace{2mm}

By Fact 1, Claims 4 and 5, Theorem 3.5 holds.
\end{proof}

\section{\normalsize Further results}\setcounter{equation}{0}

In this section, we shall study the spectral moments of graphs with given chromatic number, this parameter has closely relationship with clique number of graphs. Let $\mathscr{M}_{n, \chi}$ be the set of all $n$-vertex connected graphs with chromatic number $\chi$. Note that $\mathscr{M}_{n,n}=\{K_n\}$, hence we only consider $2\leq \chi<n$.

\begin{lem} [\cite{D-M3}]
Suppose the chromatic number $\chi(G)=t\geq4,$ let $G$ be a $t$-critical graph on more than $t$ vertices (so $G\neq K_t$). Then
$|E_G|\geq(\frac{t-1}{2}+\frac{t-3}{2(t^2-2t-1)})|V_G|.$
\end{lem}
\begin{lem}
For any $G \in\mathscr{M}_{n,t}\,(4\leq t< n),$ then $|E_G| \geq
\frac{t(t-1)}{2}+n-t,$ the equality holds if and only if $G$ is an $n$-vertex graph which is obtained from $K_t$ by attaching some trees to $K_t.$
\end{lem}

\begin{proof}
In order to determine the lower bound on the size of $G$ in $\mathscr{M}_{n,t}\,(4\leq t< n),$ it suffices to consider that $G$ is obtained from a $t$-critical graph $G'$ by attaching some trees to it. If $G'\cong K_t$, by direct computing our result holds; otherwise, consider the function
$$
   f(x)=(\frac{t-1}{2}+\frac{t-3}{2(t^2-2t-1)})x+n-x,
$$
where $t$ is a fixed positive integer with $4\le t< x$. It is easy to see that
$$f'(x)=\frac{t-1}{2}+\frac{t-3}{2(t^2-2t-1)}-1=\frac{t(t-3)(t-2)}{2(t^2-2t-1)}>0$$
for $t\geq4$. Hence, $f(x)$ is a strict increasing function in $x$, where $4\le t<x<n$. Together with Lemma 4.1 we have
\begin{eqnarray*}
|E_G|& = &|E_{G'}|+(n-|V_{G'}|)\\
     &\geq & \left(\frac{t-1}{2}+\frac{t-3}{2(t^2-2t-1)}\right)|V_{G'}|+(n-|V_{G'}|)\\
     &>& \left(\frac{t-1}{2}+\frac{t-3}{2(t^2-2t-1)}\right)t+(n-t)\\
     &>&\frac{t(t-1)}{2}+(n-t).
\end{eqnarray*}
This completes the proof.
\end{proof}

In order to determine the first few graphs in $\mathscr{M}_{n,2}$, by Lemma 2.2 these graphs must be $n$-vertex trees. Note that Pan et al. \cite{X-F-P} identified the first $\sum_{k=1}^{\lfloor\frac{n-1}{3}\rfloor}\left(\lfloor\frac{n-k-1}{2}\rfloor-k+1\right)+1$ graphs, in the $S$-order, of all trees with $n$ vertices; these $\sum_{k=1}^{\lfloor\frac{n-1}{3}\rfloor}\left(\lfloor\frac{n-k-1}{2}\rfloor-k+1\right)+1$ trees are also the first $\sum_{k=1}^{\lfloor\frac{n-1}{3}\rfloor}\left(\lfloor\frac{n-k-1}{2}\rfloor-k+1\right)+1$ graphs in the $S$-order among  $\mathscr{M}_{n,2}$. We will not repeat it here.

For convenience, let $C_n^t:=C_nu\cdot vP_{t+1}$, where $u\in V_{C_n}$, $v$ is an end-vertex of path $P_{t+1}$.
\begin{thm}
Among the set of graphs $\mathscr{M}_{n,3}$ with $n\ge 5,$ the first two graphs in the $S$-order are $C_n,\,C_{n-2}^2$ if $n$ is odd and $C_{n-1}^1,\,C_{n-3}^3$ otherwise.
\end{thm}
\begin{proof}
In order to determine the first two graphs in $\mathscr{M}_{n,3}$ with $n\ge 5$, based on Lemma 2.2 it suffices to consider the $n$-vertex
connected graphs each of which contains a unique odd cycle. Denote the set of such graphs by $\mathscr{U}_n$.

Choose $G\in \mathscr{U}_n$ such that it is as small as possible according to the $S$-order.  On the one hand, $G$ contains a unique odd cycle, say $C_t$; that is to say, $G$ is obtained by planting some trees to $C_t$ if $t<n$. On the other hand, in view of Lemma 2.5, $G$ is obtained from $C_t$ by attaching a path $P_{n-t+1}$ to it, i.e., $G\in \{C_t^{n-t}: t=3,5,\ldots\}.$

If $n$ is odd, then it suffices for us to compare $C_n$ with $C_t^{n-t}$, where $3\leq t\leq n-2$. In fact,
$S_i(C_t^{n-t})-S_i(C_n)=0$ for $i=0,\,1,\,2,\,3$. By Lemma 2.3(i), we have
$$
  S_4(C_t^{n-t})-S_4(C_n)=4[\phi_{C_t^{n-t}}(P_3)-\phi_{C_n}(P_3)]=4(n+1-n)=4>0.
$$
Hence, $C_n\prec_sC_t^{n-t}.$ Furthermore, for any $C_t^{n-t}$ with $3\leq t\leq n-4$, it is routine to check that
$S_i(C_t^{n-t})=S_i(C_{n-2}^2)=0$ for $i=0,\,1,\,2,\,3,\,4$. By direct computing (based on Lemma 2.3), we have
$$
\begin{array}{llll}
  S_5(C_t^{n-t})-S_5(C_{n-2}^2)>0 & if\ t=3, n\ge 7;& S_5(C_t^{n-t})-S_5(C_{n-2}^2)>0 & if\ t=5, n\geq9; \\
  S_5(C_t^{n-t})-S_5(C_{n-2}^2)=0 & if\ t\geq7,n\geq11; & S_6(C_t^{n-t})-S_6(C_{n-2}^2)=0 & if\ t\geq7, n\geq11; \\
  S_7(C_t^{n-t})-S_7(C_{n-2}^2)>0 & if\ t =7, n\geq11; & S_7(C_t^{n-t})-S_7(C_{n-2}^2)=0 & if\ t\geq 9, n\geq13; \\
  S_8(C_t^{n-t})-S_8(C_{n-2}^2)>0 & if\ t \geq 9, n\geq13. &  &
\end{array}
$$
This gives $C_{n-2}^2\prec_sC_t^{n-t}$. Therefore, $C_n,\,C_{n-2}^2$ are the first two graphs in the $S$-order among $\mathscr{M}_{n,3}$ for odd $n$.

If $n$ is even, it suffices to consider the graphs $C_t^{n-t}$ with $3\leq t\leq n-1$. In fact, for any $C_t^{n-t}$ with $3\leq t\leq n-3$,
$S_i(C_t^{n-t})=S_i(C_{n-1}^1)$ for $i=0,\,1,\,2,\,3,\,4$. By Lemma 2.3, we have
$$
\begin{array}{llll}
  S_5(C_t^{n-t})-S_5(C_{n-1}^1)>0 & if\ t=3, n\ge 6;& S_5(C_t^{n-t})-S_5(C_{n-1}^1)>0 & if\ t=5, n\geq 8; \\
  S_5(C_t^{n-t})-S_5(C_{n-1}^1)=0 & if\ t\geq7,n\geq10; & S_6(C_t^{n-t})-S_6(C_{n-1}^1)>0 & if\ t\geq7, n\geq10.
\end{array}
$$
Hence, we have $C_{n-1}^1\prec_sC_t^{n-t}$ for $3\leq t\leq n-1$.

Next we compare $C_{n-3}^3$ with $C_t^{n-t}$, where $3\leq t\leq n-5$. Obviously,
$S_i(C_t^{n-t})=S_i(C_{n-3}^3)$ for $i=0,\,1,\,2,\,3,\,4$. By Lemma 2.3, we have
$$
\begin{array}{llll}
  S_5(C_t^{n-t})-S_5(C_{n-3}^3)>0, & if\ t=3, n\ge 8;& S_5(C_t^{n-t})-S_5(C_{n-3}^3)>0, & if\ t=5, n\geq 10; \\
  S_5(C_t^{n-t})-S_5(C_{n-3}^3)=0, & if\ t\geq7,n\geq12; & S_6(C_t^{n-t})-S_6(C_{n-3}^3)=0, & if\ t\geq7, n\geq12;\\
  S_7(C_t^{n-t})-S_7(C_{n-3}^3)>0, & if\ t=7,n\geq12; & S_7(C_t^{n-t})-S_7(C_{n-3}^3)=0, & if\ t\geq9,n\geq14;\\
  S_8(C_t^{n-t})-S_8(C_{n-3}^3)=0, & if\ t\geq9,n\geq14. &  &
\end{array}
$$

Hence, in what follows we need compare $S_9(C_{n-3}^3)$ with that of $S_9(C_t^{n-t})$ for $t\geq9, n\geq14$.
Note that by Proposition 1, we have
$\mathscr{A}'_9(C_t^{n-t})=\{C_9\},\,\mathscr{A}'_9(C_{n-3}^3)=\emptyset,\,\mathscr{A}_9(C_t^{n-t})=\mathscr{A}_9(C_{n-3}^3)=\emptyset$
if $t=9,\,n\geq14$. Hence,
$S_9(C_t^{n-t})-S_9(C_{n-3}^3)=18(\phi_{C_t^{n-t}}(C_9)-0)=18>0$
if $t=9,\,n\geq14$, i.e., $C_{n-3}^3 \prec_s C_t^{n-t}$ in this case.

If $t\geq11,\,n\geq16$, then $\mathscr{A}'_9(C_t^{n-t})=\mathscr{A}'_9(C_{n-3}^3)=\emptyset,\,\mathscr{A}_9(C_t^{n-t})=\mathscr{A}_9(C_{n-3}^3)=\emptyset.$
Hence, $S_9(C_t^{n-t})=S_9(C_{n-3}^3)$ for $t\geq11,\,n\geq16$. Note that
$\mathscr{A}'_{10}(C_t^{n-t})=\mathscr{A}'_{10}(C_{n-3}^3)=\emptyset$,
$|\mathscr{A}_{10}(C_t^{n-t})|-|\mathscr{A}_{10}(C_{n-3}^3)|=\phi_{C_t^{n-t}}(P_6)-\phi_{C_{n-3}^3}(P_6)=n+4-(n+3)=1$
for $t\geq11,\,n\geq16$, hence
$S_{10}(C_t^{n-t})-S_{10}(C_{n-3}^3)=10(\phi_{C_t^{n-t}}(P_6)-\phi_{C_{n-3}^3}(P_6))=10>0,$ which implies  $C_{n-3}^3\prec_sC_t^{n-t}$ for $3\leq t\leq n-5$.

This completes the proof.
\end{proof}

Note that for Tur\'an graph $T_{n,t}$, $\chi(T_{n,t})=t$ and its size attains the maximum among $\mathscr{M}_{n,t}$. Combining with Lemmas 2.2 and 4.2, we have  \begin{thm}
\begin{wst}
\item[{\rm (i)}] For any graph $G\in \mathscr{M}_{n,t}\setminus\{K_t^{n-t}\}$ with $4\leq t< n$, one has $K_t^{n-t}\prec_sG.$
\item[{\rm (ii)}] For any graph $G\in \mathscr{M}_{n,t}\setminus\{T_{n,t}\}$, where $2\leq t < n$, one has $G\prec_sT_{n,t}.$
\end{wst}
\end{thm}

By Theorem 4.3, the Tur\'an graph $T_{n,t}$ is the last graph in the $S$-order among $\mathscr{M}_{n,t}$, in view of Lemma~2.2 and Theorem~3.3, the next result follows immediately.
\begin{thm}Among the set of graphs $\mathscr{M}_{n,t}$ with $2\le t \le n-1.$
\begin{wst}
\item[{\rm (i)}]If $n=t+1$, then all graphs in the set $\mathscr{M}_{n,n-1}$ have the following $S$-order:
$T_{n-3}\prec_s T_{n-4}\prec_s \cdots \prec_s T_i \prec_s \cdots\prec_s T_2\prec_sT_1\prec_s T_{n,n-1}.$
\item[{\rm (ii)}]If $n=kt$ with $3\leq t\leq\frac{n}{2}$, then for all $G\in \mathscr{M}_{n,t}\backslash\{T_{n,t},\,T_{n,t}^1\}$ one has $G\prec_s T_{n,t}^1\prec_s T_{n,t}.$
\item[{\rm (iii)}]If $n=kt+1$ with $3\leq t\leq\frac{n}{2}$, then for all $G\in \mathscr{M}_{n,t}\backslash\{T_{n,t},\,T_{n,t}^1,\,T_{n,t}^2\}$ one has $G\prec_s T_{n,t}^1\prec_s T_{n,t}^2\prec_s T_{n,t}.$
\item[{\rm (iv)}]If $n=kt+r$ with $3\leq t\leq\frac{n}{2},\, r=t-1$ or $\frac{n+1}{2}\leq t\leq n-2$, then for all $G\in \mathscr{M}_{n,t}\backslash\{T_{n,t},\,T_{n,t}^2,\,T_{n,t}^3\}$ one has $G \prec_s T_{n,t}^2 \prec_s T_{n,t}^3\prec_s T_{n,t}.$
\item[{\rm (v)}]If $n=kt+r$ with $4\leq t\leq\frac{n}{2},\, 2\leq r\leq t-2$, then for all $G\in \mathscr{M}_{n,t}\backslash\{T_{n,t},\,T_{n,t}^1,\,T_{n,t}^2,\,T_{n,t}^3\}$ one has $G\prec_s T_{n,t}^1\prec_s T_{n,t}^2 \prec_s T_{n,t}^3\prec_s T_{n,t}.$
\end{wst}
\end{thm}

By Theorems 3.4-3.5 and Lemma 4.2, we have
\begin{thm}
\begin{wst}
\item[{\rm (i)}] For $t=n-2\ge 4$, the first three graphs in the $S$-order in the set $\mathscr{M}_{n,n-2}$ are as follows:
$J_1\prec_s J_2 \prec_s J_3.$
\item[{\rm (ii)}] For $t=n-3\ge 4$, the first seven graphs in the $S$-order among the set of graphs $\mathscr{M}_{n,n-3}$ are as follows:
$J_4\prec_s J_5 \prec_s J_6 \prec_s J_7 \prec_s J_8 \prec_s J_9 \prec_s J_{10}.$
\item[{\rm (iii)}] For $4\leq t\leq n-4,$ the first $1+\sum_{i=1}^{\lfloor\frac{n-t-1}{3}\rfloor}(n-t-3i)$ graphs, in the $S$-order, among the set of graphs $\mathscr{M}_{n,t}$ are as follows:
$K_t^{n-t}\prec_sG^1_{n-t-2}\prec_sG^1_{n-t-3}\prec_s\cdots\prec_sG^1_3\prec_sG^1_{2}\prec_sG^2_{n-t-4}\prec_sG^2_{n-t-5}\prec_s\cdots\prec_sG^2_3\prec_s\cdots
\prec_s G^i_{n-t-2i} \prec_s G^i_{n-t-2i-1}\prec_s \cdots \prec_s G_{i+1}^i \prec_s G^i_{n-t-2(i+1)} \prec_s \cdots \prec_s
G^{\lfloor\frac{n-t-1}{3}\rfloor}_{\lfloor\frac{n-t-1}{3}\rfloor+1}.$
\end{wst}
\end{thm}

\end{document}